\documentclass{amsart}

\usepackage{amssymb}
\usepackage[all]{xy}
\usepackage{hyperref}
\newtheorem{thm}{Theorem}%[section]

\newtheorem{lem}[thm]{Lemma}

\newtheorem{conj}[thm]{Conjecture}
\theoremstyle{definition}

\newtheorem{say}[thm]{}

\newtheorem*{ack}{Acknowledgments}      % \renewcommand{\theack}{} 

\newtheorem{defn-thm}[thm]{Definition--Theorem}  %!!!!!!!!!!!!!!!!!!!!!!!!
\newtheorem{defn-lem}[thm]{Definition--Lemma}  %!!!!!!!!!!!!!!!!!!!!!!!!
\theoremstyle{remark}

\newcommand{\z}[0]{{\mathbb Z}}

\newcommand{\p}[0]{{\mathbb P}}

\newcommand{\q}[0]{{\mathbb Q}}

\newcommand{\qtq}[1]{\quad\mbox{#1}\quad}

\newcommand{\supp}[0]{\operatorname{Supp}}

\newcommand{\sing}[0]{\operatorname{Sing}}

\newcommand{\simq}[0]{\sim_{\q}}

\def\into{\DOTSB\lhook\joinrel\to}

\def\loccoh#1.#2.#3.#4.{H^{#1}_{#2}(#3,#4)}

\DeclareMathAlphabet{\mathchanc}{OT1}{pzc}%
                                {m}{it}

\usepackage[all]{xy}\xyoption{dvips}

\begin{document}
\bibliographystyle{amsalpha}
\hfill\today

 \title{Moduli of polarized Calabi-Yau pairs}
 \author{J\'anos Koll\'ar and Chenyang Xu}

 \maketitle

Compactifying the moduli space of Calabi-Yau varieties is a challenging problem. 
For a family of K3 surfaces over a punctured disc, Kulikov \cite{MR0506296}  discovered that, after a base change,  there are 
degenerations with trivial canonical class, but these are usually reducible and  there are   infinitely many non-isomorphic ones.

The same general framework   holds in higher dimensions as well. 
There are degenerations with trivial canonical class and semi-log-canonical singularities, but usually infinitely many non-isomorphic ones. 
It does not seem possible to choose one in a sensible and functorial way, unless one imposes extra structures.
In increasingly  general forms these claims were proved in
\cite{ksb, ale-pairs, bchm, MR3032329, NX16, k-nx}.

A general approach to obtain unique degenerations was first explored by Alexeev \cite{ale-abvar} for Abelian varieties and Hacking \cite{hacking} for plane curves. Instead of  Calabi-Yau varieties one needs to  work with pairs  $(X, \epsilon H)$, where $X$ is a 
Calabi-Yau variety,  $H$ is an ample divisor on $X$ and $0<\epsilon\ll 1$. Note that $H$ is a divisor, not a linear equivalence class or a cohomology class. However, in many cases there is a distinguished  choice of the divisor $H$, and then 
this approach is especially natural and useful. This happens in the papers
\cite{ale-abvar, hacking, 2019arXiv190309230D, 2019arXiv190309742A}.

 For maximal  generality, we  
 consider 
 {\it polarized Calabi-Yau pairs.} These consist of a   {\it Calabi-Yau pair} (that is, a semi-log-canonical pair  $(X, \Delta)$, as in  \cite[Chap.5]{kk-singbook}, where 
$K_X+\Delta$ is $\q$-linearly equivalent to 0), plus 
an ample $\q$-Cartier divisor $H$ on $X$ such that $(X, \Delta+\epsilon H)$  is also semi-log-canonical for
$0<\epsilon\ll 1$. 

Polarized Calabi-Yau pairs have a natural moduli space, denoted by
$\mathbf{PCY}$,  see Paragraph~\ref{pcy.defn}. 
By  construction  $\mathbf{PCY}$ is locally of finite type,
but it has infinitely many connected components since we did not even fix the dimension of $X$.  Thus the best one could hope for in general is the following.

\begin{conj} The connected components of $\mathbf{PCY}$ are projective.
\end{conj}

Note that this somewhat goes against the conjectures of Reid \cite{MR909231},
but the two are not inconsistent since the polarized deformations considered here are not the same as the non-polarized deformations studied in \cite{MR909231}.

The aim of this note is to prove a weaker statement, which is however usually sufficient in all concrete situations.
This theorem seems to have been known to several people, but it was listed as an open question in some recent preprints, for example in \cite{2019arXiv190309230D, 2019arXiv190309742A}.   Thus it  may be worthwhile to write  down the precise statement of the general result and its proof. 
(In \cite{2019arXiv190309230D, 2019arXiv190309742A} the emphasis is on describing  particular irreducible components of  $\mathbf{PCY}$  in concrete terms, so none of the main results of these papers are effected by our theorem.)

\begin{thm}\label{mpcy.thm}
 The irreducible components of $\mathbf{PCY}$ are projective.
\end{thm}

First we need to  fix the definitions of the relevant moduli problems. From now on we work over a field $k$ of characteristic 0.

\begin{say}[Moduli of stable pairs] (See \cite{k-modsurv, k-modbook} for details.)

A {\it pair} $(X,\Delta)$ consist of a reduced, pure dimensional variety $X$ and an effective  $\q$-divisor $\Delta$, none of whose irreducible components is contained in $\sing X$. For moduli purposes it is best to  write
$\Delta=c D$, where $D$ is a $\z$-divisor. 

A pair $(X,\Delta)$ is called  {\it locally stable} if it is 
semi-log-canonical, see \cite[Chap.5]{kk-singbook}.  $(X,\Delta)$ is {\it stable } if, in addition,  $X$ is projective and  $K_X+\Delta$ is an ample $\q$-divisor.

The general definition of  {\it locally stable} and  {\it stable}  morphisms is somewhat complicated. However, if  $S$ is normal, then the following works, see \cite[Chap.3]{k-modbook}. 

A morphism  $f:(X,\Delta)\to S$ is {\it locally stable} if $f:X\to S$ is flat, $K_{X/S}+\Delta$ is $\q$-Cartier and  all fibers are locally stable.
(Implicitly,  this includes the condition that the fibers should make sense, that is, the restriction of $\Delta$ to any fiber should make sense.
The later holds iff $\Delta$ is $\q$-Cartier at the generic points of
$X_s\cap\supp\Delta$ for every $s\in S$. There are several technical issues with this when $S$ is not normal, but these become crucial only when we pass from one irreducible component of the moduli space to another. Thus these
 are not relevant for our current purposes.  See  \cite[Chap.4]{k-modbook} for a discussion.)

A morphism   is {\it  stable}  if, in addition, 
  $f$ is  projective and 
$K_{X/S}+\Delta$ is $f$-ample.

The main theorem is that, at least in characteristic 0,
 there is a coarse moduli space
$\mathbf{SP}$, which is separated and satisfies the valuative criterion of properness.

Once the existence of  $\mathbf{SP}$ is shown, 
the boundedness results of \cite{al94, MR3779687}  imply that 
the connected components of $\mathbf{SP}$ are proper, and then
\cite{k-proj, MR3779955, MR3671934} imply that they are projective. 
\end{say}

\begin{say}[Moduli of polarized Calabi-Yau pairs]\label{pcy.defn}
  Working in the most general setting,
a {\it Calabi-Yau pair} is a proper, semi-log-canonical pair  $(X, \Delta)$ where 
$K_X+\Delta$ is $\q$-linearly equivalent to 0. We write
$\Delta=c D$, where $D$ is a $\z$-divisor and consider the constant $c$ as fixed in our moduli problem.

 For example, following \cite{hacking}, when we work with the moduli of hypersurfaces $D$ of degree $d\geq n+1$ in $\p^n$, then we think of the objects as Calabi-Yau  pairs
$$
\bigl(X, \Delta=\tfrac{n+1}{d}\cdot D\bigr),
$$
where $D$ is allowed to be reducible and can even have components with multiplicity $\geq 2$ for $d\geq 2(n+1)$. 

A {\it polarized Calabi-Yau pair} consists of a Calabi-Yau pair $(X, \Delta)$
plus an ample $\q$-Cartier divisor $H$ such that 
$(X, \Delta+\epsilon H)$ is semi-log-canonical for $0<\epsilon\ll 1$. 
The latter holds iff $H$ does not contain any of the semi-log-canonical centers of $(X, \Delta)$, see \cite[2.5 and 2.13]{kk-singbook}.

A {\it stable family} of polarized Calabi-Yau pairs over a normal 
base scheme $S$ consists of
  a flat, proper morphism  $f:X\to S$, 
 a $\q$-divisor  $\Delta$ on $X$ and
  a $\q$-Cartier divisor $H$
 such that
$K_{X/S}+\Delta$ is $\q$-Cartier and all fibers 
$(X_s, \Delta_s, H_s)$ are polarized Calabi-Yau pairs.

Let us now fix a rational $0<\epsilon< 1$ and consider those  polarized Calabi-Yau pairs for which $(X, \Delta+\epsilon H)$ is semi-log-canonical.
Then  $(X, \Delta+\epsilon H)$ is a  stable pair. The corresponding objects
form an open subset 
$$
\mathbf{PCY}_{\epsilon}\subset \mathbf{SP},
$$
that gives the moduli space of those polarized Calabi-Yau pairs for which 
$(X, \Delta+\epsilon H)$ is semi-log-canonical.

If we pick a smaller $0<\epsilon_2<\epsilon_1 $ then the sets
$\mathbf{PCY}_{\epsilon_1}$ and $\mathbf{PCY}_{\epsilon_2}$ are actually disjoint  (since we decreased the self-intersection of
$K_X+\Delta+\epsilon H$), but 
sending $(X,\Delta,\epsilon_1 H)$ to $(X,\Delta,\epsilon_2 H)$ defines an 
 open embedding
$$
j(\epsilon_1, \epsilon_2): \mathbf{PCY}_{\epsilon_1}\into \mathbf{PCY}_{\epsilon_2}.
$$
As $\epsilon\to 0$, the directed union of these embeddings 
gives the {\it moduli space of polarized Calabi-Yau pairs.} We denote it by
$\mathbf{PCY}$. 
\end{say}

\begin{say}[Difficulties of the traditional approach to Theorem~\ref{mpcy.thm}]
Assume for simplicity that $\Delta=0$.

Working with one Calabi-Yau variety $X$, we   take  $(X, \epsilon H)$ with $\epsilon$ small enough, its precise value is not important.
However, the value of $\epsilon$ becomes crucial in families.

Consider a family $f:(X, \epsilon H)\to S$ in $\mathbf{SP}$ whose generic fiber is in 
$\mathbf{PCY}_{\epsilon}$. That is, $f:X\to S$ is  a flat, projective morphism,
 $K_{X/S}+\epsilon H$ is $f$-ample and $K_{X/S}$ is trivial on the  generic fiber.  Three problems can happen if we want to change $\epsilon$.
\begin{itemize}
\item If $H$ is not $\q$-Cartier then any change in $\epsilon$ results in a family that is not allowed in our moduli theory.
% \item Increasing $\epsilon$ may results in a family $f:(X,\eta H)\to S$
% whose fibers  $(X_s, \eta H_s)$ are not log canonical.
\item Decreasing $\epsilon$ may result in a family $f:(X,\eta H)\to S$
for which $K_{X/S}+\eta H$ is not $f$-ample.
\end{itemize}
While both are known to happen for some values of $\epsilon$, standard conjectures of 
the theory of minimal models suggest that if we start with any  family
$f:(X,\epsilon H)\to S$ and gradually decrease the value of  $\epsilon$, then, 
after finitely many contractions and  flips we should get a new family
$f^m:(X^m,\epsilon_0 H^m)\to S$ where neither of the above problems occur for any further decrease of $\epsilon_0$. 

This should give  a very satisfactory answer for any given  family, but there is one more problem.
\begin{itemize}
\item The value of $\epsilon_0$ may need to  get arbitrarily small, depending on the family we start with, even for families with the same generic fiber.
 \end{itemize}
The latter  is usually referred to as a boundedness question of the corresponding moduli problem. In our case the 
  general  boundedness results of \cite{al94, MR3779687}  do not apply since the 
underlying varieties are Calabi-Yau and  the value of $\epsilon$ is not fixed. 
\end{say}

We solve the first 2 problems by first  running a carefully chosen auxiliary MMP as in \cite{MR3032329}. Then we note that any 
irreducible component of $\mathbf{PCY}$ is covered by a single universal family, so we evade the third problem as well.

\begin{say}[Proof of Theorem~\ref{mpcy.thm}]
We prove that the irreducible components of $\mathbf{PCY}$ are proper. Then the general results of
\cite{k-proj, MR3779955, MR3671934} imply that they are projective.

 Let $\mathbf M$ be an irreducible component of $\mathbf{PCY}$
with generic point $g_M$. Then there is a finite extension of
$K\supset k(g_M)$ such that we have a 
polarized Calabi-Yau pair $(X_K, \Delta_K, H_K)$ over $K$ that corresponds to $g_M$. We prove in Lemma~\ref{lem1} that for  $0<\epsilon\ll 1$
there is a projective variety $S$ such that $k(S)$ is a finite extension of $K$, and  a stable family of
polarized Calabi-Yau pairs
$$
f_{S}:\bigl(X_{S}, \Delta_{S}+\epsilon H_{S}\bigr)\to S,
$$
such that over the generic point we recover $(X_{k(S)}, \Delta_{k(S)}, H_{k(S)})$. %possiby after a further finite extension of $K$.

If this holds then consider the moduli map $S\to \mathbf{PCY}$.
Its image contains $g_M$ and it is proper 
since $S$ is projective. Thus $\mathbf M$, which is the closure of $g_M$ in $\mathbf{PCY}$, is proper. \qed
\end{say}

\begin{lem}\label{lem1} Let  $K/k$ be a function field and $(X_K, \Delta_K, H_K)$ a polarized Calabi-Yau pair over $K$. Then there is a projective variety $S$ such that $k(S)/K$ is finite  and  a stable family of
polarized Calabi-Yau pairs
$$
f_{S}:\bigl(X_{S}, \Delta_{S}+\epsilon H_{S}\bigr)\to S
$$
extending $(X_K, \Delta_K+ \epsilon H_K)\times_Kk(S)$ for some for $0<\epsilon\ll 1$.
\end{lem}

Proof. 
Assume first that  $X_K$ is normal  and geometrically irreducible.
Choose  a log resolution  $\pi_K:(Y_K, \Delta^Y_K+H^Y_K)\to (X_K, \Delta_K+H_K)$
such that 
$$
K_{Y_K}+ \Delta^Y_K\simq \pi_K^*\bigl(K_{X_K}+ \Delta_K\bigr)\qtq{and}
H^Y_K=\pi_K^*H_K.
$$
We can extend it to 
 a simultaneous   log resolution
$$
(Y_{S_1},\Delta^Y_{S_1}+H^Y_{S_1})\to (X_{S_1}, \Delta_{S_1}+H_{S_1})
$$
over some affine variety $S_1$ such that $k(S_1)\cong K$.   By  \cite[Thm.0.3 and Sec.8.2]{Abramovich-Karu00}  there is a projective, generically finite, dominant morphism $\pi:S_2\to S_1$
and a compactification  $S_2\into S$ such that
the pull-back $\bigl(Y_{S_1},\supp(\Delta^Y_{S_1}+H^Y_{S_1})\bigr)\times_{S_1}S_2$ extends to a locally stable morphism
$$
g_{S}:\bigl(Y_{S}, \supp(\Delta^Y_{S}+H^Y_{S})\bigr)\to S,$$ where  $S$ is smooth
and $Y_{S} $ has only quotient (hence  $\q$-factorial) singularities. 
Note that every log canonical center of $\bigl(Y_{S}, \supp(\Delta^Y_{S}+H^Y_{S})\bigr)$ dominates $S$.

Due to the presence of the  quotient  singularities, we can not guarantee that 
$\bigl(Y_{S}, \supp(\Delta^Y_{S}+H^Y_{S})\bigr)$ be dlt. However, it has {\it qdlt singularities} (that is quotients of dlt singularities) as discussed in \cite[Sec.5]{dkx}. 

Write $\Delta^Y_{S}=\Theta^+_S-\Theta^-_S$ as the difference of effective divisors without  common irreducible components.
By construction  $(X_{k(S)}, \Delta_{k(S)})$ is a good minimal model of
 $(Y_{k(S)}, \Theta^+_{k(S)}) $.

Let $ L^Y_S$ be a general, sufficiently ample divisor on $Y_S$. Then, for
$0<\eta\ll 1$, 
\begin{enumerate}
\item $L^Y_S+H^Y_S$ is relatively ample, 
\item $\bigl(Y_{S}, \Theta^+_{S}+\eta(L^Y_S+H^Y_S)\bigr)$ is qdlt and 
\item  $g_S: \bigl(Y_{S}, \Theta^+_{S}+\eta(L^Y_S+H^Y_S)\bigr)\to S$ locally stable.
\end{enumerate}
By  \cite[1.1]{MR3032329} and \cite[2.9]{MR3779687} 
the relative minimal model program with scaling of $L^Y_S+H^Y_S$  for
   $g_{S}:(Y_{S}, \Theta^+_{S})\to S$
terminates with 
  $$
g^m_{S}:\bigl(Y^m_{S}, \Delta^m_{S}+\epsilon(L^m_S+H^m_S)\bigr)\to S
$$
for some $0<\epsilon\leq \eta$ 
such that 
$K_{Y^m_{S}/S}+\Delta^m_{S}$ is $g^m_{S}$-semiample.
(Note that $K_{Y_{S}/S}+ \Theta^+_{S} $ is $\q$-linearly equivalent to 
$\Theta^-_{S} $ on the generic fiber. Thus the above MMP contracts $\supp \Theta^-_{S} $, so the images of $\Theta^+_{S} $ and of $\Delta_S^Y$ agree on $Y^m_S$. This is why we  change notation back to $\Delta^m_{S} $.)

A relative minimal model of a projective, locally stable  morphism might not be locally stable, but this holds if the base space is smooth by
\cite[Cor.10]{k-nx}. 
Thus $$
g^m_{S}:\bigl(Y^m_{S}, \Delta^m_{S}+\epsilon (L^m_S+H^m_S)\bigr)\to S
\qtq{ is locally stable.}
$$  
The log canonical class $K_{Y^m_{S}/S}+\Delta^m_{S}$  is semiample
and of Kodaira dimension 0 on  the generic fiber, hence it is relatively $\q$-linearly trivial.
Thus $g^m_{S}:(Y^m_{S}, \Delta^m_{S})\to S$
 is a locally stable family of Calabi-Yau pairs.
Since  $Y_{S} $ is $\q$-factorial, so is  $Y^m_{S} $. In particular,
$H^m_S$ is $\q$-Cartier and so
 $$
g^m_{S}:\bigl(Y^m_{S}, \Delta^m_{S}+\epsilon H^m_S\bigr)\to S
\qtq{ is also locally stable.}
$$  

% From $H_K$ by pull-back we get $\pi_K^*H_K$. After base extension we get a divisor on the generic fiber of $Y_S$; let  $H^Y_S$ denote its closure.
% By assumption $(X_K, \Delta_K+ \epsilon H_K)$ is log canonical 
% for $0<\epsilon\ll 1$, so the pair $(Y_{S}, \Delta^Y_{S}+\epsilon H^Y_S) $
% is also log canonical, even dlt,  along the generic fiber. As we noted, 
%  every log canonical center of $(Y_{S}, \Delta^Y_{S}) $ dominates $S$,
% thus in fact 
% $(Y_{S}, \Delta^Y_{S}+\epsilon H^Y_S) $ is  qdlt for $0<\epsilon\ll 1$.

% Let now $H^m_S$ denote the birational transform of $H^Y_S$ on $Y^m_{S}$.
% Any given step of a $(Y_{S}, \Delta^Y_{S}) $-MMP is also a
% step of a   $(Y_{S}, \Delta^Y_{S}+\epsilon H^Y_S) $-MMP  for $0<\epsilon\ll 1$,
% thus $(Y^m_{S}, \Delta^m_{S}+\epsilon  H^m_S)$ is qdlt  for $0<\epsilon\ll 1$. 

Since $H^m_S$ is relatively big,
by \cite[1.1]{MR3032329}  $(Y^m_{S}, \Delta^m_{S}+\epsilon  H^m_S)\to S$ has a relative canonical model
$$
f_{S}:\bigl(X_{S}, \Delta_{S}+\epsilon H_{S}\bigr)\to S.
$$
As before, \cite[Cor.10]{k-nx} guarantees that 
$f_{S}$ is stable. 
Since $K_{Y^m_{S}/S}+\Delta^m_{S}$ is relatively $\q$-linearly trivial,
the same holds for $K_{X_{S}/S}+\Delta_{S}$. Thus
$$
f_{S}:\bigl(X_{S}, \Delta_{S}+\epsilon H_{S}\bigr)\to S
$$
is a  family of
polarized Calabi-Yau pairs. This completes the proof when
$X_K$ is normal.

If $X_K$ is not normal, let  $(X^n_{K}, \Delta^n_{K}+\epsilon H^n_K) $ denote its normalization. 
The previous step, applied to each irreducible component,  gives  $f^n_{S}:(X^n_{S}, \Delta^n_{S}+\epsilon H^n_S)\to S$. Finally the gluing theory of \cite{k-source}  and \cite[Chap.5]{kk-singbook}
 applies and  we get 
$f_{S}:(X_{S}, \Delta_{S}+\epsilon H_S)\to S$. See \cite[Sec.2.4]{k-modbook} for details. 
\qed

\begin{ack} 
We thank V.~Alexeev  and K.~DeVleming for helpful conversations. 
Partial  financial support    was provided  by  the NSF under grant numbers
 DMS-1362960 to JK and  DMS-1901849  to CX. Both authors also received support  from the grant  DMS-1440140 while  in residence at
MSRI during the Spring 2019 semester.
\end{ack}

%\bibliography{refs}

\def\cprime{$'$} \def\cprime{$'$} \def\cprime{$'$} \def\cprime{$'$}
  \def\cprime{$'$} \def\dbar{\leavevmode\hbox to 0pt{\hskip.2ex
  \accent"16\hss}d} \def\cprime{$'$} \def\cprime{$'$}
  \def\polhk#1{\setbox0=\hbox{#1}{\ooalign{\hidewidth
  \lower1.5ex\hbox{`}\hidewidth\crcr\unhbox0}}} \def\cprime{$'$}
  \def\cprime{$'$} \def\cprime{$'$} \def\cprime{$'$}
  \def\polhk#1{\setbox0=\hbox{#1}{\ooalign{\hidewidth
  \lower1.5ex\hbox{`}\hidewidth\crcr\unhbox0}}} \def\cdprime{$''$}
  \def\cprime{$'$} \def\cprime{$'$} \def\cprime{$'$} \def\cprime{$'$}
\providecommand{\bysame}{\leavevmode\hbox to3em{\hrulefill}\thinspace}
\providecommand{\MR}{\relax\ifhmode\unskip\space\fi MR }
\providecommand{\MRhref}[2]{%
  \href{http://www.ams.org/mathscinet-getitem?mr=#1}{#2}
}
\providecommand{\href}[2]{#2}

\bigskip

\noindent  Princeton University, Princeton NJ 08544-1000

\email{kollar@math.princeton.edu}
\medskip

\noindent   MIT, Cambridge, MA 02139-4307

\email{cyxu@mit.edu}

\end{document}